\documentclass[12pt, a4paper]{article}

\usepackage{amsmath, amsthm, amscd, amsfonts, amssymb, graphicx, color}
\usepackage[bookmarksnumbered, plainpages, backref]{hyperref}
\usepackage{graphicx}
\usepackage{graphicx,epstopdf}
\usepackage{tikz}
\usepackage{float}
\usepackage[all]{xy}
\newcommand{\chain}[1][n]{\langle #1\rangle}
\newcommand{\zetaf}[1][P]{\mathfrak{Z}\left(#1\right)}

\newcommand{\nsomega}[2][n]{\Omega(#2,#1)}
\newcommand{\somega}[2][n]{\Omega^{\circ}(#2,#1)}

\newcommand{\bast}{\hat *}
\newcommand{\dd}{\widehat D}
\newcommand{\cupdot}{\mathbin{\mathaccent\cdot\sqcup}}
\def\multiset#1#2{\left(\kern-.3em\binom{#1}{#2}\kern-.3em\right)}

\usepackage{hyperref}
\newcommand{\footremember}[2]{%
    \footnote{#2}
    \newcounter{#1}
    \setcounter{#1}{\value{footnote}}%
}
 
\title{Operad of Posets 101: The Wixárika Posets}

\author{%
  José Antonio Arciniega-Nevárez\footremember{alley}{División de Ingenierías, Universidad de Guanajuato, Guanajuato, México. \vskip 0.1cm }%
  \and Marko Berghoff\footremember{trailer}{Research Data Centre, German Federal Pension Insurance, Berlin, Germany. \vskip 0.1cm}%
  \and Eric Rubiel Dolores-Cuenca\footremember{al}{ Yonsei University, Seoul, Korea.}%
  }

\begin{document}

\maketitle

\begin{abstract}

We study classes of objects whose combinatorics are closely related to those of posets. The framework of operads and operad algebras allows us to make this relationship precise and provides tools for a deeper understanding of their combinatorial structure. 

In this note, we present a nontrivial example of a suboperad of the operad of posets, called Wix\'arika posets, together with its associated algebras. This example is sufficiently rich to exhibit key structural features of the theory, while remaining accessible and avoiding unnecessary technicalities.

\end{abstract}

\newtheorem{theorem}{Theorem}[section]
\newtheorem{lemma}[theorem]{Lemma}
\newtheorem{proposition}[theorem]{Proposition}
\newtheorem{corollary}[theorem]{Corollary}
\theoremstyle{definition}
\newtheorem{definition}[theorem]{Definition}
\newtheorem{example}[theorem]{Example}
\newtheorem{xca}[theorem]{Exercise}
\newtheorem{problem}[theorem]{Problem}
\theoremstyle{remark}
\newtheorem{remark}[theorem]{Remark}
\numberwithin{equation}{section}

\tableofcontents
\section{Introduction}

This note is intended as an extension of lecture notes designed to assist students in understanding the theory of operads, with a particular focus on algebras over the operad of posets. This pedagogical aim is inspired by the well-known role of Abdus Salam as both a teacher and a promoter of science, emphasizing the importance of making advanced mathematical structures accessible to students and early-career researchers. Accordingly, this note does not contain new results.

Let us begin with an overview of what this note is about. Let $P$ be a partially ordered set (poset). We associate to $P$ a directed acyclic simple graph, the \emph{Hasse diagram} $H_P$ of $P$: Each vertex of $H_P$ represents an element of $P$ and there is an arrow from $x$ to $y$ if  $x< y$ and there is no $z$ with $ x< z< y$. We can visualize the composition of posets $P$ and $Q$ in terms of their Hasse diagrams. Given a point  $p\in P$, the composition $P\circ_p Q$ is the poset whose Hasse diagram is the Hasse diagram of $P$ with the vertex $p$ replaced by the diagram $H_Q$. 

Using this composition, we can associate to the poset $\{x,y\}$, the disjoint union of two elements, a map of posets that sends $P,Q$ to their disjoint union,
\begin{equation*}
\{x,y\}_{\hbox{Poset}}(P,Q)=(\{x,y\}_{\hbox{Poset}}\circ_x P)\circ_y Q=P\sqcup Q.
\end{equation*}

This algebraic structure can be formalized using the language of operads.
An operad consists of sets $O(n)$ for all $n\geq 1$, together with composition maps $O(n)\times O(k_1)\times\cdots\times O(k_n)\mapsto O(\sum_{i=1}^n k_i)$, satisfying certain compatibility conditions (see Section~\ref{Sec:op}). In the previous example, the operad of finite posets is formed by $\mathrm{Poset}(n)$, the sets of posets with $n$ elements. For $P\in \mathrm{Poset}(n), Q_i\in \mathrm{Poset}(k_i), i=1,\cdots,n,$ the composition $P\circ (
Q_1,\cdots,Q_n
)$ is the poset whose Hasse diagram is obtained by replacing for each $i=1,\ldots,n$ the $i$-th vertex of $H_P$ by $H_{Q_i}$. 

Thus, from the operadic point of view each poset also represents an endomorphism of posets. Motivated by this construction, we aim to study sets $X$ together with maps $\mathrm{Poset}(n)\mapsto \mathrm{Map}(X^n,X)$. For a correspondence that sends a poset $P$ to a map $f_P$, we also require that the poset $P\circ_{p_i} Q$ is sent to the map composition $f_P(id,\cdots,id, f_Q,id,\cdots,id)$ where $f_Q$ is located on the $i$-th entry.  

Here is an explicit example: R.\ Stanley \cite{beginning} defined the \emph{strict order polynomial} $\Omega^\circ(P,n)$ of $P$ by counting the number of strict order-preserving maps from $P$ to the ordered sets $\chain[n]=\{1<\cdots<n\}, n \in \mathbb{N}$. By definition, $\Omega^\circ(P\sqcup Q, n)=\Omega^\circ(P, n)\Omega^\circ(Q, n)$. We rephrase this result by saying that the endomorphism associated with $\{x,y\}$ 
is the product of strict order polynomials, that is, 
\begin{equation*}
    \{x,y\}_{\hbox{Polynomials}}(\Omega^\circ(P, n),\Omega^\circ( Q, n))=\Omega^\circ
\left(\{x,y\}_{\hbox{Poset}}(P, Q), n \right) =\Omega^\circ(P, n)\Omega^\circ(Q, n).
\end{equation*}
Note that here and in what follows we include a subscript to indicate the set on which a poset acts.

Another example comes from polytope theory. The order polytope of a poset $P$ is defined as 
\begin{equation*}
\mathrm{Poly}(P) = \{ f \colon P \longrightarrow [0,1] \mid x \leq_P y \Rightarrow f(x)\leq f(y) \}.   
\end{equation*}
One can verify that 
\begin{align*}
\mathrm{Poly}(\{x,y\}_{\hbox{Poset}}(P,Q)) & =\{x,y\}_{\hbox{Polytope}}(\mathrm{Poly}(P),\mathrm{Poly}(Q)) \\
&=\{ v+w \mid v \in \mathrm{Poly}(P), w \in \mathrm{Poly}(Q) \},
\end{align*}
that is, the map $\{x,y\}_{\hbox{Polytope}}$ is the Minkowski sum of polytopes. 
 
From the operadic point of view, we study sets in which the elements of the operad define maps. These sets are called algebras over the operad. Above we described order polynomials and order polytopes as examples of algebras over the operad of posets. It will become clear that posets themselves are an algebra over the operad of posets, and we will learn how to transfer information from one algebra to the other. 

As this is an introductory text, it will be more convenient to work with a particular family of posets which we call Wixárica posets. We outline the theory for this case. As an application, we describe how topological information of these posets (i.e., of their Hasse diagrams) is reflected at the algebraic level.

The exposition is organized as follows. In Section~\ref{Sec:2} we define posets and introduce Wix\'arika posets. In Section~\ref{sec:order-series} we introduce an algebra over the operad of Wix\'arika posets: Generating functions of the above mentioned order polynomials, called \emph{order series}. Section~\ref{Sec:op} formalizes the notion of operads and their algebras.  
Section~\ref{sec:algo} provides a worked-out example in enumerative combinatorics. Section~\ref{sec:repr} gives an example of how topological information is transferred from the Wix\'arika operad to one of its algebras. We then suggest future reading in Section~\ref{Section:final}.

\section{Wixárica Posets}\label{Sec:2}
In this section, we define our main object of study, the class of Wixárica posets. We are interested in algebraic properties of these posets. In particular, we will define the operations of concatenation and disjoint union of posets. Moreover, we will introduce a new operation, called the 'itsari operation, to obtain Wixárica posets through the consecutive application of the latter operations. We then interpret these in terms of Hasse diagrams of Wixárica posets in order to obtain graphical representations of these operations.

In principle, posets may be defined using two variants of order relations, either strict or non-strict orders. Most of the time, we work with the latter. We shall give both definitions for completeness.  

\begin{definition}
A \emph{partial order} $\leq$ on a set $P$ is a reflexive, antisymmetric and transitive relation on its elements. That is, for $a,b,c\in P$:
\begin{enumerate}
\item \emph{Reflexivity}: $a\leq a$, i.e., every element is related to itself.
\item \emph{Antisymmetry}: if $a\leq b$ and $b\leq a$, then $a=b$, i.e., no two distinct elements precede each other.
\item \emph{Transitivity}: if $a\leq b$ and $b\leq c$, then $a\leq c$. 
\end{enumerate}
\end{definition}

\begin{definition}
A \emph{strict partial order} $<$ on a set $P$ is an irreflexive, asymmetric and transitive relation on its elements. That is, for $a,b,c\in P$:
\begin{enumerate}
\item \emph{Irreflexivity}: no element $a$ is related to itself.
\item \emph{Asymmetry}: if $a<b$, then it is not possible that $b<a$.
\item \emph{Transitivity}: if $a< b$ and $b< c$, then $a< c$. 
\end{enumerate}
\end{definition}

\begin{definition}
A set $P$ together with a (strict) partial order $(P,\leq)$ is called a \emph{(strict) partially ordered set} or simply a \emph{(strict) poset}. We use $P$ to denote both the underlying set and the poset. If confusion is possible, we add a subscript (e.g., $\leq_P$) to distinguish different partial orders. 
\end{definition}

\begin{example}Examples of posets:
 Ordered sets, such as the real numbers; the set of subsets of a given set, its power set, ordered by inclusion; the set of natural numbers equipped with the divisibility relation; if we consider a finite  set of programs with the order relations $A<B$ whenever the output of $A$ is an input of $B$, then a set of programs forms a poset where the order gives restrictions to select at what turn to run each program.
\end{example}

In this note we assume that every poset is finite.

A \emph{Hasse diagram} is a graph that represents a finite poset by drawing its transitive reduction. We will draw Hasse diagrams when we need to visualize certain posets. Concretely,  a partially ordered set $P$ with order $\leq$  is represented by a graph whose vertices are the elements of $P$, and we have an edge from one vertex $r$ to another vertex $s$ whenever $r\neq s$ and $r< s$ and there is no $t$ with $r< t< s$. Typically, in Hasse diagrams the greater element is positioned higher than the lesser. 

\begin{figure}[H]
    \centering
    \includegraphics[width=0.5\linewidth]{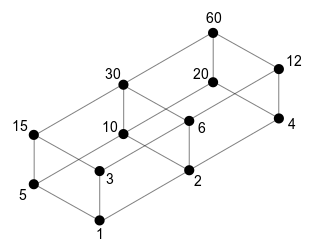}
    \caption{Hasse diagram of the poset of all divisors of 60. Edges denotes the order relation, and higher position encodes the hierarchy (Image created by Ed\_g2s, \url{https://commons.wikimedia.org/wiki/User:Ed_g2s}).}
    \label{fig:divisor}
\end{figure}

\begin{example}
An example is given in Figure~\ref{fig:divisor}. It shows the poset given by the divisors of 60. There is an edge from 1 to the numbers 2, 3, and 5, since 1 is their closest divisor. For the remaining numbers, we can apply transitivity. This means that 1 is a divisor of 5, and 5 is a divisor of 15, but an edge between 1 and 15 is not required. Note also that 1 is located at the bottom, since it is the smallest number. Each number is placed on a higher level according to its order. For this reason, 60 is located at the top.
\end{example}

We use a particular class of posets, \emph{series-parallel posets}. These posets are defined constructively, starting from a set of points and iterated applications of operations that we now introduce.

\begin{definition}\label{Def*}
We define the \emph{direct sum} of two posets $(P_1,\leq_{P_1})$ and $(P_2,\leq_{P_2})$ to be the poset $P=P_1*P_2$ whose underlying set is the union of $P_1\sqcup P_2$ and $a\leq_P b$ if and only if 
\begin{enumerate}
\item $a,b\in P_1$ and $a\leq_{P_1} b$, or
\item $a,b\in P_2$ and $a\leq_{P_2} b$, or
\item $a\in P_1$ and $b\in P_2$.
\end{enumerate}
\end{definition}

\begin{definition}
The disjoint union poset of $(P_1,\leq_{P_1})$ and $(P_2,\leq_{P_2})$ is the poset whose underlying set is $P_1\sqcup P_2$ and there are no more relations than those given by $\leq_{P_1}$ and $\leq_{P_2}$.
\end{definition}

Recall from the introduction how the disjoint union $P\sqcup Q$ can be formulated using poset composition and the additional poset $\{x,y\}$, a disjoint union of two elements:  
\begin{equation*}
P\sqcup Q = \{x,y\}_{\hbox{Poset}}(P,Q)=(\{x,y\}_{\hbox{Poset}}\circ_x P)\circ_y Q.
\end{equation*}
In the same way the poset $\{x<y\}$ gives rise to the direct sum,
\begin{equation*}
    P * Q = \{x < y\}_{\hbox{Poset}}(P,Q)=(\{x < y\}_{\hbox{Poset}}\circ_x P)\circ_y Q.
\end{equation*}

From now on we stick to the symbols $\sqcup$ and $*$ to keep the notation light. Nevertheless, it is instructive in what follows to keep in mind that these operations on posets are induced by posets themselves.

In the context of series-parallel posets, the direct sum is called \emph{series composition} while \emph{parallel composition} refers to the disjoint union.

\begin{definition}
The class of \emph{series-parallel orders} is defined inductively as follows:
\begin{enumerate}
\item A single point is series-parallel.
\item If $(P_1,\leq_{P_1})$ and $(P_2,\leq_{P_2})$ are series parallel with $P_1\cap P_2=\emptyset$, then their series
composition is series-parallel.
\item If $(P_1,\leq_{P_1})$ and $(P_2,\leq_{P_2})$ are series-parallel with $P_1\cap P_2=\emptyset$, then their parallel composition is series-parallel. 
\end{enumerate}
\end{definition}

\begin{figure}[H]
\centering
\includegraphics[width=0.5\linewidth]{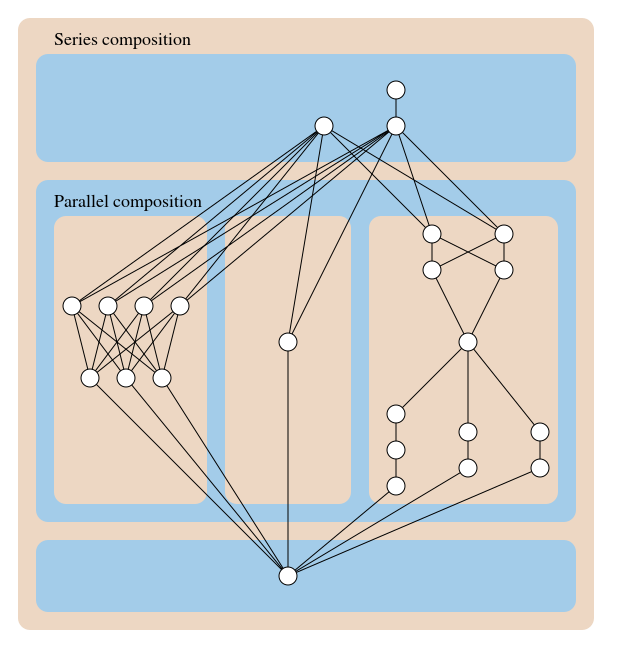}
\caption{Hasse diagram of a series-parallel poset. Posets in brown rectangles are joined together by disjoint union (parallel composition), posets in blue rectangles are joined by direct sum (series composition). (Image created by David Eppstein, \url{https://commons.wikimedia.org/wiki/User:David_Eppstein/Gallery}).}
\label{fig:series-parallel}
\end{figure}

\begin{example}
An example is illustrated in Figure~\ref{fig:series-parallel}. The central blue rectangle illustrates parallel composition, that is, a disjoint union of three posets, one in each brown rectangle. The final poset is obtained by two series compositions, where in each case edges are added from the maximal element(s) of one blue rectangle to the minima element(s) of the next blue rectangle.     
\end{example}

Next, we introduce a new operation on single posets. Consecutive applications of this operation and disjoint union will produce the Wixárica posets. Note that Wixárica posets are a special case of series-parallel posets. 

\begin{definition}\label{ioperation}
The \emph{'itsari} operation $D$ is defined  as
\begin{equation*}
D(P)=\chain[1]*\big(\chain[1]\sqcup P\big)*\chain[1].
\end{equation*}
In the case of chains $\chain$ we abbreviate $D(\chain)$ by $D(n)$ (cf.\ Section \ref{sec:algo}).
\end{definition}
'Itsari is a word from Wixárica language that means to weave. Geometrically, we are adding a new minimum $x_0$ and a new maximum $x_1$ to $X$, then we attach a handle on these points, that is, we add a third element $y$ with $x_0<y<x_1$.

\begin{definition}
A \emph{Wixárika poset} is a finite poset obtained from the poset with one element by finitely many applications of the operations $D$ and $*$.    
\end{definition}

\begin{example}
An example of a Wixárica poset is shown by the Hasse diagram in Figure~\ref{fig:melted}. We start by building a chain using the direct sum. Then we add arcs using the 'itsari operation. Note that in contrast to the previous examples, we direct Hasse diagrams of Wixárica posets from left to right.
\end{example}
\begin{figure}[H]
\centering
\begin{tikzpicture}[scale=0.7]
\coordinate  (v1) at (-7,0);  
\coordinate  (v2) at (-6,0);
\coordinate  (v3) at (-5,0);
\coordinate  (v4) at (-4,0);
\coordinate  (v5) at (-3,0);
\coordinate  (v6) at (-2,0);
\coordinate  (v7) at (-1,0);
\coordinate  (v8) at (0,0);
\coordinate  (v9) at (1,0);
\coordinate  (v10) at (2,0);
\coordinate  (v11) at (3,0);
\coordinate  (v12) at (4,0);
\coordinate  (v13) at (4.5,0);
\coordinate  (v14) at (5,0);
\coordinate  (v15) at (5.5,0);
\coordinate  (v16) at (6,0);
\coordinate (w1) at (-3,-.7);
 \coordinate  (w2) at (1,-.7);  
 \coordinate  (x1) at (-1,-2);
  \coordinate  (y1) at (-1,-3);
  \coordinate  (z1) at (-1,-4);
  
      \draw (v1) -- (v16);
\draw (v4) to[out=-90, in=180] (w1); \draw (w1) to[out=0, in=-90] (v6);
\draw (v8) to[out=-90, in=180] (w2); \draw (w2) to[out=0, in=-90] (v10);
\draw (v3) to[out=-90, in=180] (x1); \draw (x1) to[out=0, in=-90] (v11);
\draw (v2) to[out=-90, in=180] (y1); \draw (y1) to[out=0, in=-90] (v12);
\draw (v1) to[out=-90, in=180] (z1); \draw (z1) to[out=0, in=-95] (v16);

\fill[black] (v1) circle (.0666cm);
\fill[black] (v2) circle (.0666cm);
\fill[black] (v3) circle (.0666cm);
\fill[black] (v4) circle (.0666cm);
\fill[black] (v5) circle (.0666cm);
\fill[black] (v6) circle (.0666cm);
\fill[black] (v7) circle (.0666cm);
\fill[black] (v8) circle (.0666cm);
\fill[black] (v9) circle (.0666cm);
\fill[black] (v10) circle (.0666cm);
\fill[black] (v11) circle (.0666cm);
\fill[black] (v12) circle (.0666cm);
\fill[black] (v13) circle (.0666cm);
\fill[black] (v14) circle (.0666cm);
\fill[black] (v15) circle (.0666cm);
\fill[black] (v16) circle (.0666cm);
\fill[black] (w1) circle (.0666cm);
\fill[black] (w2) circle (.0666cm);
\fill[black] (x1) circle (.0666cm);
\fill[black] (y1) circle (.0666cm);
\fill[black] (z1) circle (.0666cm);
\end{tikzpicture}
\caption{Hasse diagram of a Wixárika poset, oriented from left to right.}
\label{fig:melted}
\end{figure}

The Hasse diagrams of these posets resemble Wix\'arika collars, which are intricate colorful necklaces crafted by the Wixárika community in North America.

\section{Order Series}\label{sec:order-series}
In this section, we describe a set of power series associated to posets and operations between them. The idea is that these operations are related to the poset operations defined in the previous section. 

Most of our examples of posets and ordered sets use a simple type of poset, called a chain.

\begin{definition}
Define a \emph{$n$-chain} to be the total ordered set $\chain=\{1< 2<\cdots <n\}$. Let $P$ be a finite partially ordered set. An \emph{order-preserving map} $\sigma:P\to \chain$ is a map such
that if $r\leq s$ in $P$, then $\sigma(r)\leq\sigma(s)$. Define $\Omega(P,n)$ as the number of order-preserving maps from $P$ to $\chain$. A \emph{strictly order-preserving map} $\tau:P\to\chain$ is a map such that if $r<s$ then $\tau(r)<\tau(s)$. Define $\Omega^\circ(P,n)$ to be the number of strictly order-preserving maps from $P$ to $\chain$.
\end{definition}
\begin{definition}Given $n,m\in\mathbb{N}$ define the multiset $
    \multiset{n}{m}:={n+m-1\choose m}$.
\end{definition}
\begin{example}
If $P=\chain[m]$ we have $\somega[n]{\chain[m]}={n\choose m}$, while $\nsomega[n]{\chain[m]}=\multiset{n}{m}$. If $P$ is the disjoint union of $p$ points, then $\somega[n]{P}=\nsomega[n]{P}=n^p$.
\end{example}

\begin{definition}
For a poset $P$ we define a formal power series, called the \emph{(strict) order series} of $P$, by
\begin{equation*}
\zetaf[P]=\zetaf[P,x]=\sum_{n= 1}^\infty \somega[n]{P}x^n.
\end{equation*}
\end{definition}

\begin{example}
For instance, we have
\begin{gather}
\zetaf[m]:=\zetaf[{\chain[m]}]=
 \sum_{n=m}^\infty {n\choose m} x^n =\frac{x^m}{(1-x)^{m+1}}.
  \label{eqn:other1}
\end{gather} 
For a proof of the second equality see \cite[Equation~(1.5.5)]{gf}, \cite[Equation~(1.3)]{lcomb} or \cite[Equation~(1.3)]{eulerianB}.
\end{example}

We will use some standard operations on power series with order series. We include them for completeness.

\begin{definition}
If $\sum_{i=0}^\infty a_ix^i$ and $\sum_{i=0}^\infty b_ix^i$ are two power series, the \emph{Cauchy product} is defined as
\begin{equation*}
\left(\sum_{i=0}^\infty a_ix^i\right)\left(\sum_{i=0}^\infty b_ix^i\right)=\sum_{n=0}^\infty c_nx^n\quad \text{where} \quad c_n=\sum_{k=0}^na_kb_{n-k}.
\end{equation*}
\end{definition}

\begin{definition}\label{def:astseries}
If $x$ is the common variable for the order series $\zetaf[P]$ and $\zetaf[Q]$, we define the \emph{direct product} using the Cauchy product as follows: 
\begin{equation}
\zetaf[P]\bast\zetaf[Q]=\zetaf[P](1-x)\zetaf[Q].\label{eq:*}  
\end{equation}
\end{definition}

\begin{definition}
We denote the \emph{Hadamard product} of power series by
\begin{equation*}\sum_{n=1}^\infty a_n x^n \odot \sum_{n=1}^\infty b_n x^n=\sum_{n=1}^\infty a_n b_n x^n.
\end{equation*}     
In the case of order series, we use the symbol $\cupdot$, that is,
\begin{equation*}
\zetaf[P]\cupdot\zetaf[Q]=\zetaf[P]\odot\zetaf[Q].
\end{equation*}
\end{definition}

\begin{definition}
We define the 'itsari operation for power series as
\begin{equation*}
\dd (\zetaf[P])=\zetaf[D(P)].
\end{equation*}
In the case of chains $\chain$ we abbreviate $\dd(\zetaf[n])$ by $\dd(n)$ (cf.\ Section \ref{sec:algo}).
\end{definition}

The following two formulas are well-known (see \cite{algo}, \cite{prior},  
\cite{refl},
\cite{survey},\cite{power}).

\begin{proposition}\label{prop:os-concatenation}
If $P$ and $Q$ are two posets, then
\begin{gather}
\zetaf[P*Q]=\zetaf[P]\bast\zetaf[Q]=\zetaf(1-x)\zetaf[Q],\label{eq:direct-produc}\\
\zetaf[P\sqcup Q]=\zetaf[P]\cupdot\zetaf[Q]=\zetaf[P]\odot\zetaf[Q]\label{eq:direct-sum}.
\end{gather}
\end{proposition}

Note that Equation~\eqref{eq:direct-sum} follows directly from counting the number of order preserving maps from $P\sqcup Q$ to $n$. In the case of chains, the paper~\cite{power} provides an explicit formula for the order series of disjoint unions:
 
\begin{proposition}\label{prop:os-disjoin-u}
For two chains we have
\begin{equation*}
\zetaf[{\chain[k]} \sqcup {\chain[m]}]=\sum_{n=0}^k{m+n \choose k}{k \choose n}\zetaf[m+n]. 
\end{equation*}
\end{proposition}
\begin{corollary}\label{cor:zetacadenas}
The strict order series of any series-parallel poset is a linear combination of order series of chains.\label{Cor:lc}
\end{corollary}

The following result will be very useful for one of our main objectives, which is to describe the order series of any series-parallel poset.

\begin{proposition}\label{Prop:linear}
Unraveling the definitions, we can show linearity of $\dd$: \begin{equation*}
\dd(a\zetaf[P]+b\zetaf[Q])=a\dd(\zetaf[P])+b\dd(\zetaf[Q]).
\end{equation*}
\end{proposition} 

\begin{example}\label{ex:d-operad}
If the input of the 'itsari operation $D$ is an $n$-chain, we compute
\begin{align*}
\zetaf[D(n)]&=\zetaf[1]\bast\zetaf[{\chain[1]\sqcup \chain[n]}]\bast\zetaf[1]
\\
&=\frac{x}{(1-x)^2}\bast\left(n\frac{x^n}{(1-x)^{n+1}}+(n+1)\frac{x^{n+1}}{(1-x)^{n+2}}\right)\bast\frac{x}{(1-x)^2}\\
&=n\frac{x^{n+2}}{(1-x)^{n+3}}+(n+1)\frac{x^{n+3}}{(1-x)^{n+4}}\\
&=n\zetaf[n+2]+(n+1)\zetaf[n+3].
\end{align*}
\end{example}

\section{Algebras over the operad of Wix\'arika posets}\label{Sec:op}
In this section, we introduce the language of operads. This allows us to compare sets in which the same operations are defined, similarly to groups and group homomorphisms.

\begin{definition}
An operad $(\mathcal{O}, \circ )$ consists of a collection $\mathcal{O}=\{O(n)\}_{n\geq 1}$ of sets and an associative composition $\circ : O(n) \times O(k_1)\times\cdots\times O(k_n)\to O(\sum_{i=1}^n k_i)$. There are unital elements $id\in O(1)$, and there is an action of $S_n$ that is compatible with the composition operation. 
\end{definition}
In the composition, an element in $O(n)$ takes the role of a map with $n$ inputs, and an element in $O(k_1)\times\cdots\times O(k_n)$ represents $n$ maps that are composed with the former. The result is a map with $\sum_{i=1}^n k_i$ inputs . 
\begin{definition}
    
Given an operad $\mathcal{O}$ and $f\in O(n)$, for every $i\leq n$ and $g\in O(k_i)$ the partial composition $\circ_i$ is defined by $f\circ_i g= f\circ(id,\cdots,id,g,id,\cdots,id)$  where $g$ is in the position of the $i$-th input of $f$. 
\end{definition}

A more detailed definition of operad can be found in~ \cite{whatis,OpTre,ope, alghom, whatare,entropy,dendroidal}.   

\begin{example}\label{ex:end-operad}
The most fundamental example of an operad is
the \emph{endomorphism operad} $\mathrm{End}_A$: Let $A$ be a set and ${\mathrm{Map}(A^n,A)}$ the set of functions from the cartesian product $A^n$ to $A$. Then $\mathrm{End}_A$ has as $n$-ary elements $\mathrm{End}_A(n)={\mathrm{Map}(A^ n,A)}$.
\end{example}

\begin{example}Consider the \emph{tree grafting operad}: Let $T(n)$ be the set of trees with $1$ root and $n$ leaves labeled (arbitrarily) $1$ through $n$. The collection  $\{T (n)\}$  forms an operad by grafting the root of $g$ to the leaf of $f$
labeled $i$, as in Figure~\ref{fig:f-operad}. Once the grafting is done, the leaves of $g$ are relabeled adding $i-1$ to the original label. Leaves of $f$ after the leaf $i$ are also relabeled continuing the numeration $i+m,\cdots,m+n-1$.
\begin{figure}[H]
\centering
\includegraphics[width=\linewidth]{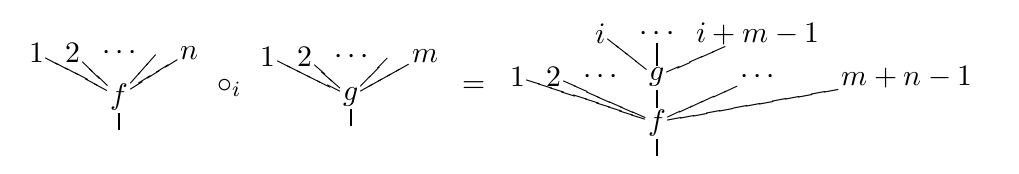}
\caption{Operations on the operad tree. The tree $g$ is grafted in the tree $f$ using the leaf positioned at $i$ (Image created by Wiki\_cies, \url{https://commons.wikimedia.org/w/index.php?title=User:Wiki_cies&action=edit&redlink=1}).}
\label{fig:f-operad}
\end{figure}
\end{example}

\begin{definition}
    By an algebra $A$ over an operad $\mathcal{O}$ we mean an operadic morphism $\mathcal{O}\mapsto \mathrm{End}_A$.
\end{definition}

An algebra $A$ over an operad $\mathcal{O}$ is then a set in which the abstract elements of $O(n)$ are realized as functions $A^n\mapsto A$. Heuristically, we can think of elements $\sigma\in O(n)$ of an operad as $n$-ary trees. Then a composition $\sigma\circ_i\tau$ corresponds to the grafting of the trees were we identify the leaf of $\sigma$ with the root of $\tau$.  The evaluation $\sigma(a_1,\cdots, a_n)$ where  $a_1,\cdots,a_n$ are elements of $A$ then corresponds to placing drops of water  on the leaves of the tree. The drops fall along the tree. Whenever a drop reaches a branching point (labeled with an operation), we evaluate the corresponding operation and replace the drop by the output. The process continues until we reach the root. 

\begin{definition}
Let $(P,\leq_P)$ be a finite poset with cardinality $|P|=n$. Let $\mathbf{Q}=\{(Q_i,\leq_{Q_i})\}_{i\in P}$ be a collection of finite posets indexed by the elements of $P$. The lexicographic sum of elements in the collection $\mathbf{Q}$ along $P$ is the poset with underlying set $\bigcup Q_i$ and the order $x\leq y$ is given by:
\begin{enumerate}
\item  if $x,y\in Q_i$ and  $x\leq_{Q_i}y$ or
\item  if $x\in Q_i$, $y\in Q_j$, and  $i\leq_{P}j$.
\end{enumerate}
\end{definition}

\begin{example} Recall that we defined the binary endomorphism $\ast$ on the set of finite posets in Definition~\ref{Def*}. Given finite posets $Q_x$ and $Q_y$, the lexicographic sum of $Q_x$ and $Q_y$ along $\{x<y\}$ is $Q_x\ast Q_y$.
\end{example}

\begin{definition}
The \emph{Wixárica operad} is the operad $W$ generated under lexicographic sum by the unary operation $D$, the identity, and the binary operation $\ast$.      
\end{definition}

\begin{example}
The elements of the  operad $W$ can be described by binary rooted trees with bird nests: There are vertices with two incoming edges and one outgoing edge labeled by $*$; and the bird nest stand for vertices with one incoming and one outgoing edge, which  are labeled by $D$. See Figure \ref{fig:tree}.

\begin{figure}[H]
\centering{
\begin{tikzpicture}[scale=.6]
\coordinate  (v1) at (-5,1);  
\coordinate  (v2) at (-5,-1);
      \draw (v1) -- (v2);
      
\coordinate  (w1) at (-3,1);
\coordinate  (w2) at (-3,-1);
\coordinate  (w3) at (-3,0);   
\draw (w1) -- (w2);
\fill[black] (w3) circle (.18cm);

\coordinate  (x1) at (-1.6,1);
\coordinate  (x2) at (-.4,1);
\coordinate  (x3) at (-1,0);
\coordinate  (x4) at (-1,-1);
     \draw (x1) -- (x3);
      \draw (x2) -- (x3);
      \draw (x3) -- (x4);
\end{tikzpicture}
\hspace{1.5cm}
 \begin{tikzpicture}[scale=.35]
 \draw (0,0) -- (0,8);
 \filldraw[black] (0,1) circle (.18cm);
  \filldraw[black] (0,3.3) circle (.18cm);
   \filldraw[black] (0,5) circle (.18cm);
   
    \draw (0,6) -- (-1.5,7.5);
   \filldraw[black] (-.75,6.75) circle (.18cm);
   
       \draw (0,7) -- (1.5,8);
   \filldraw[black] (.75,7.5) circle (.18cm);
   
      \draw (0,1.8) -- (-2,2.8);
      \draw (-1.4,2.5) -- (-1.4,3);
      \draw (-1,2.3) -- (-1,3);
\end{tikzpicture}}
\caption{Left: Trees illustrating the identity operation, the operation $D$ and the operation $\ast$. Right: A tree representing the element in Figure~\ref{fig:melted}, which is obtained using only the operations $D$ and $\ast$.}\label{fig:tree}
\end{figure}
\end{example}

The lexicographic sum then allows us to see Wix\'arika posets as algebras over the operad  $W$ generated by points.

\begin{example}
Consider Figure~\ref{fig:des}. On the left side we have at the top three copies of the posets $\chain[1]$. Below is a sequence of transformations using $D$ and $*$. At the bottom we have the final result, a Wixárika poset.
On the right side we have the corresponding tree of operations in which the leaves should be thought of as decorated by the poset $\chain[1]$. The vertices with two incoming edges and one outgoing edge are labeled with $*$ and the vertices with one incoming edge and one outgoing edge (the bird nest) are labeled with $D$. From top to bottom, as drops of water (the posets) located on the leaves fall due to gravity, the drops reach the branches labeled with operations. When passing through, drops are replaced by the output of the operations. This continues until we reach the root.  
\begin{figure}[H]
\centering
\includegraphics[scale=0.5]{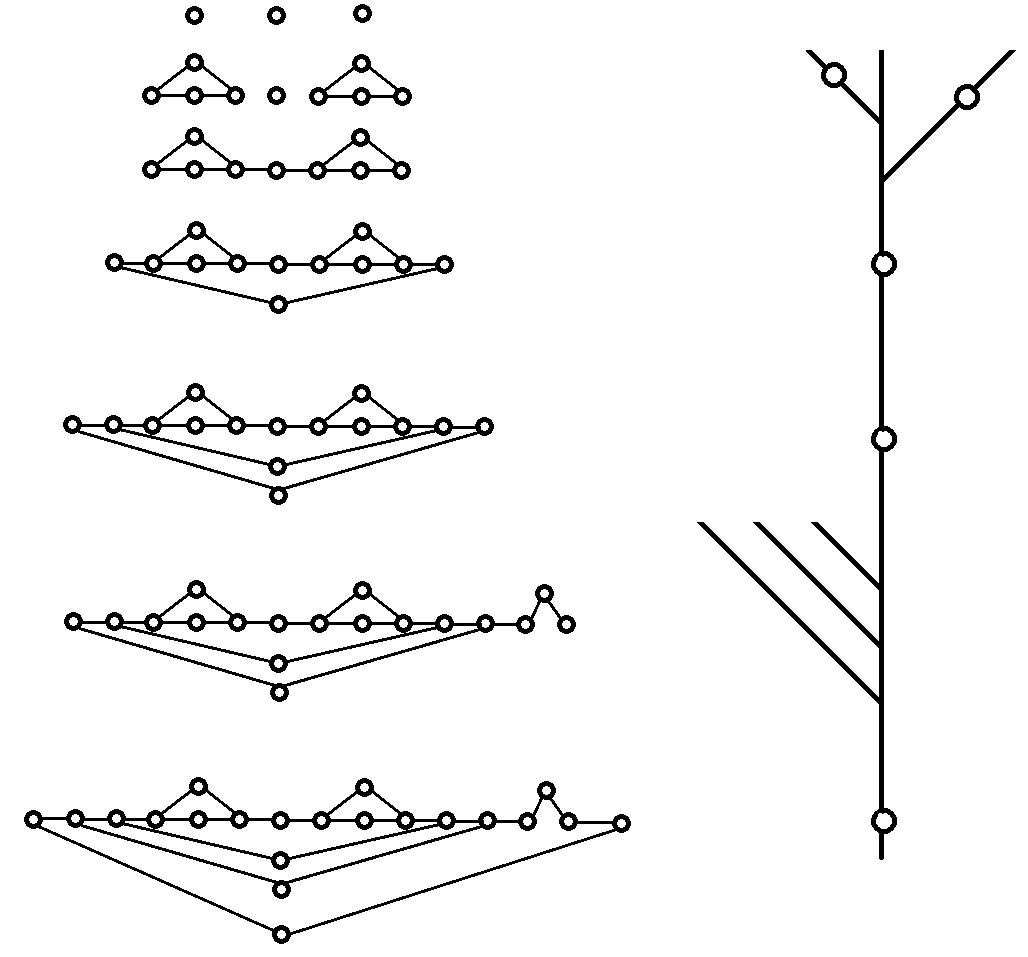}
\caption{To the left bottom we have a Wixárika poset, as an element of the algebra. To the right we have the tree of operations, an element of the operad of Wix\'arika posets. When this tree of operations is evaluated on six points it produces the Wix\'arika poset of the left. As the drops of water run down the tree, the poset on the left is being assembled.}\label{fig:des}
\end{figure}     
\end{example}

\begin{definition}
Given two algebras $A,B$ over an operad $\mathcal{O}$, an operadic morphism $\phi:A\to B$ is a map that makes the following diagram commutative:
\begin{equation*}
\xymatrix{\mathcal{O}\ar[rd]\ar[d]&\\\mathrm{End}_A\ar[r]&\mathrm{End}_B\\}
\end{equation*}
\end{definition}

\begin{proposition}\label{proposition:structure}
The function $\zetaf[\cdot]$ defines an operadic homomorphism between the $W$-algebra of Wix\'arika posets and the $W$-algebra of order series of Wix\'arika posets.     
\end{proposition}

\begin{proof}
This proposition is another way to paraphrase Equations~\eqref{eq:direct-produc} and \eqref{eq:direct-sum}. Put together Definition~\ref{def:astseries} and Proposition~\ref{prop:os-concatenation} to get:
\begin{equation*}\zetaf[P\ast Q]
=\zetaf[P](1-x)\zetaf[Q]=\zetaf[P]\bast\zetaf[Q]. \end{equation*}    
By Proposition~\ref{prop:os-disjoin-u} and its Corollary~\ref{cor:zetacadenas} we have:
\begin{equation*}
\zetaf[P\sqcup Q]=\zetaf[P]\cupdot\zetaf[Q].    
\end{equation*}
 We obtain the explicit formula
\begin{equation*}
\dd(\zetaf[P])=\zetaf[1]\bast(\zetaf[P]\cupdot \zetaf[1])\bast\zetaf[1].
\end{equation*}
\end{proof}

\section{Evaluation on order series}\label{sec:algo}
In this section, we give an explicit example of how to compute the order series of a Wixárica poset using the algebra homomorphims between both $W$-algebras defined in the previous section. To compute the order series, we follow \cite{algo}, but work on the level of power series instead of order polynomials.

Consider the homomorphism between the algebra over the  operad $W$ generated by $\chain[1]$ and the algebra generated by $\zetaf[1]$. In Proposition~\ref{proposition:structure} we constructed $\bast$ and $\cupdot$ such that
\begin{equation*}
\zetaf[P]\bast\zetaf[Q]=\zetaf[P*Q],\quad \dd(\zetaf[P])=\zetaf[D(P)].
\end{equation*}

If we have an element of the operad, which we describe as a binary tree with bird nests, then we obtain a poset by placing points $\chain[1]$ on the leaves of the tree. Similarly, we obtain a power series by placing $\zetaf[1]$ on the leaves of the tree, that is, $\zetaf[]$ substitutes any input leaf $\chain[1]$ by $\zetaf[1]$. 

\begin{example}\label{ex:op-homomorphism}
Consider the poset in Figure~\ref{fig:des}. The resulting poset is given by 
\begin{equation*}
P=D(D(D(D(1)*\chain[1]*D(1)))*\chain[3]).
\end{equation*} 
We shall give a procedure to calculate $\zetaf[P]$. The first step of our algorithm is replacement of the leaves $\chain[1]$ by their order series $\zetaf[1]=\frac{x}{(1-x)^2}$. In the tree, going from top to bottom, the algebra homomorphism allows us to make the following replacements: a concatenation $P*Q$ by $\zetaf \bast \zetaf[Q]$ and $D$ by $\dd$. Using Proposition~\ref{prop:os-concatenation} and Example~\ref{ex:d-operad}, as well as linearity of $\dd$ and bilinearity of $\bast$, one evaluates the vertex operations until we reach the root to obtain
\begin{gather*}
\dd\left(\dd\left(\dd\left(\dd(1)\bast\zetaf[1]\bast\dd(1)\right)\right)\bast\zetaf[3]\right)\\
=\dd\left(\dd\left(\dd\left(\left(\zetaf[3]+2\zetaf[4]\right)\bast\zetaf[1]\bast\left(\zetaf[3]+2\zetaf[4]\right)\right)\right)\bast\zetaf[3]\right)\\
=\dd\left(\dd\left(\dd\left(\zetaf[7]+4\zetaf[8]+4\zetaf[9]\right)\right)\bast\zetaf[3]\right)\\
=\dd\left(\dd\left(\dd(7)\right)\bast\zetaf[3]\right)+4\dd\left(\dd\left(\dd(8)\right)\bast\zetaf[3]\right)
+4\dd\left(\dd\left(\dd(9)\right)\bast\zetaf[3]\right).
\end{gather*}

Following the same procedure we continue:
\begin{equation*}
=882\frac{x^{16}}{(1-x)^{17}}+7995\frac{x^{17}}{(1-x)^{18}}+27232\frac{x^{18}}{(1-x)^{19}}+
\end{equation*}
\begin{equation*}
+43792\frac{x^{19}}{(1-x)^{20}}+33552\frac{x^{20}}{(1-x)^{21}}+9880\frac{x^{21}}{(1-x)^{22}}.\end{equation*}

As a second step, to compute the order polynomial $\Omega^\circ(m)$, we use the expression $\frac{x^n}{(1-x)^{n+1}}=\sum_{m\geq n} {m\choose n}x^m$ to extract the term of degree $x^m$ from the order series. Because our work focuses on studying operations that preserve the number of labeling maps of posets, we know that the coefficient of $x^m$ is the number of labelings of the poset using numbers from $1$ to $m$. 

From the above we thus find that for any $m$ the number of labeling maps of $P$ is given by 
\[882{m\choose 16}+7995{m\choose 17}+27232{m\choose 18}+43792{m\choose 19}+33552{m\choose 20}+9880{m\choose 21}.\]     
\end{example}

\section{Representability}\label{sec:repr}
We now use the topological information to answer the inverse problem: Given a function $f(x)$, does there exist a Wixárika poset $P$ whose order series coincides with $f(x)$? In this section, we present an algorithm to answer this question.

\begin{definition}
We say a function $f(x)$ is \emph{represented} if there is a Wixárika poset $P$ such that $f(x)=\zetaf$. 
\end{definition}

\begin{example}
Two posets are called \emph{Doppelg\"{a}nger posets} if they share the same order series, for example $\chain[1]*(\chain[1]\sqcup\chain[1]\sqcup\chain[1])$ and $\chain[2]\sqcup \chain[2]$ have the same order series: $\zetaf[2]+6\zetaf[3]+6\zetaf[4]$ \cite{Doppelgangers}.    
\end{example}

In order to decide whether a function can be represented, we determine topological invariants of a poset $P$ from its associated power series $\zetaf[P]$, such as the Betti number of the Hasse diagram of $P$ (i.e., the number of independent cycles in the underlying undirected graph). With this information, we describe an algorithm that determines whether $f(x)$ can be represented.  If this is the case, then the algorithm explicitly constructs all possible Wixárika posets $P$ such that $f(x)=\zetaf[P]$.

The following theorem states that the coefficients of the order series associated to a Wixárika poset are topological invariants of its Hasse diagram.
 
\begin{theorem}\label{prop:rep}
Let $P$ be a Wixárika poset. Consider its factorization, that is, the associated word $w$ in the characters $*$ and $D$. Let  $f(x)=w(\zetaf[1])$ be the power series obtained by replacing every leaf of $w$ by $\zetaf[1]$ and proceeding as in Example~\ref{ex:op-homomorphism}. 

\begin{enumerate}
\item The series $f(x)$ admits a unique description 
\begin{equation*}
f(x)= d_i\zetaf[i]+d_{i+1}\zetaf[i+1]+\cdots+d_{k}\zetaf[k],
\end{equation*}
where the coefficients $d_j$ are positive integers. 
\item The number $i$ is the number of points in a maximal chain in $P$.\label{prop:point}
\item The number $k$ is the number of points in $P$.
\item The difference $d=k-i$ is the first Betti number of the Hasse diagram of $P$, equal to the number of times $D$ occurs in $w$.
\item The difference $m=i-2d-1$ is the number of times  $*$ occurs in $w$. 
\item The term $m+1$ is the number of leaves in the tree $w$.
\item\label{item:cond} $\sum_{u=1}^{k}(-1)^{k-u}d_u=1$.
\end{enumerate}
\end{theorem} 
\begin{proof}[Sketch of the proof]
All statements follow by describing the possible values obtained under the operations $\dd$ and $\bast$. The explicit description of $\dd$ and $\bast$ leads to the positivity of the coefficients.

Consider, for example, statement $(4)$: Our goal is to show that the difference $k-i$ counts the first Betti number of the Hasse diagram of the poset. We can argue inductively:

The first Betti number of $\chain[1]$ is zero, and $k-i=0$ for $\zetaf[1]$.

We know that by applying $D$ we increase the first Betti number of the Hasse diagram of a poset by 1. We also know the precise effect of applying $\dd$ on $\zetaf[n]$. According to Example~\ref{ex:d-operad} the difference $k-i$ is increased by 1. By linearity (Proposition~\ref{Prop:linear}) the effect of $\dd$ on the order series of a poset increases $k-i$ by 1 as well.

Wixárika posets satisfy the following property: The first Betti number of the poset obtained by concatenation is the sum of the first Betti numbers of the input posets. This additivity also holds for the difference $k-i$: Given two power series of Wixárika posets $f,g$, from Proposition~\ref{def:astseries} we see that $f\bast g$ satisfies $(k_f+k_g)-(i_f+i_g)=(k_f-i_f)+(k_g-i_g)$. 

Finally, since Wixárika posets are generated by $\chain[1]$ under the operations $D,*$, we conclude that the difference $k-i$ counts the first Betti number of a Wixárika poset. 

For details, see \cite[Proposition 3.1]{power}.
\end{proof}
\begin{example}\label{ex:rep}
The series \begin{equation*}
882\frac{x^{16}}{(1-x)^{17}}+7995\frac{x^{17}}{(1-x)^{18}}+27232\frac{x^{18}}{(1-x)^{19}}+
\end{equation*}
\begin{equation*}
+43792\frac{x^{19}}{(1-x)^{20}}+33552\frac{x^{20}}{(1-x)^{21}}+9880\frac{x^{21}}{(1-x)^{22}},\end{equation*} corresponds to the poset depicted on the left bottom of Figure~\ref{fig:des}. We visually verify that the poset has $i=16$ points in any chain of maximal size, consists of $k=21$ elements in total, the Betti number of its Hasse diagram equals $k-i=5$, is generated by $m=16-2(5)-1=5$ uses of $*$, and its tree on the right has $6$ leaves.
\end{example}

The paper~\cite{survey} calls the elements $\{\zetaf[i]\}_{i\in\mathbb{N}}$ the \textit{chain basis}.
The coefficients of order series in the chain basis satify the following positivity property:

\begin{theorem}
\label{lemma:a}
Consider $\zetaf=\sum_{i=j_0}^{|P|} d_i\zetaf[i], d_i\geq0$. The order polytope of $P$ is the union of $d_{|P|}$ copies of the $|P|$-simplex. The remaining coefficients $d_i$ of $\zetaf=\sum d_i\zetaf[i]$ count the number of $i$-simplices that occur, in the canonical triangulation of the order polytope, as intersection of the $d_{|P|}$ copies of those $|P|$-simplices.
 \label{aux}
\end{theorem}

\begin{proof}
See \cite[Remark 2.17 and Proposition 3.6]{power}.
\end{proof}

The coefficients of power series associated to more general posets, for instance series-parallel posets (which are generated by concatenation and disjoint union), are not topological invariants. For example, the Hasse diagram of the poset $\chain[1]*(\chain[1]\sqcup \chain[1])*\chain[1]$ has first Betti number $1$ while the Hasse diagram of the poset
 $\chain[2]*(\chain[1]\sqcup \chain[1])$ has first Betti number $0$. Their order series, however, are the same.

In the other direction, we have the following statement.

\begin{theorem}\label{thm:alg}
Given a finite sum $f(x)=\sum d_j\zetaf[j]$ where $d_j$ are non-negative integers, there is an algorithm to determine if $f(x)=\zetaf$ for some Wixárika poset $P$. In the positive case the algorithm returns all possible posets representing $f(x)$.
\end{theorem}
\begin{proof}[Sketch of proof]
Given $f$ we test posets to see whether $\zetaf[P]=f$ for some poset $P$. We use geometric information to reduce the search space.

First, we verify that $f$ satisfies conditions (1) and (7) of Theorem~\ref{prop:rep}.
Using  the remaining conditions of Theorem~\ref{prop:rep}, we then filter candidate posets (e.g., Betti number equal to $k-i$, number of points equal to $i$, etc).
Next, we pass to equivalence classes of posets up to associativity (for example, by using colored operads to factorize the posets). 
We then evaluate candidate posets in stages. For every representative $P$ of an equivalence class of posets, we only evaluate the first and last element in the expansion $\zetaf[P]=\sum_{j=i}^k c_j\zetaf[j]$, that is, we keep track only of $c_i\zetaf[i]$ and $ c_k\zetaf[k]$. 

At every stage of the evaluation of $\zetaf[P]$, if we obtain partial coefficients $(c_i^\prime\zetaf[i], c_k^\prime\zetaf[k])$, we require that $c_i^\prime | d_i$, otherwise the poset can be discarded. The evaluation is performed using the tree of operations as described in \cite[Section~4]{algo}, but working with power series instead of polynomials. 

If we find a poset $Q $ with $\zetaf[Q]=r_i\zetaf[i]+\cdots+r_k\zetaf[k]$ and $r_i=d_i, r_k=d_k$ then we compute all intermediate coefficients ${\{r_j\}}_{i<j<k}$ and compare the result with $f(x)$. We continue until we find a solution or all candidate posets have been tested. If a poset producing $f(x)$ is found, we use the results of \cite[Sections 4.2 and 5 ]{survey} to determine all \textit{Doppelg\"{a}ngers}, that is, posets with the same power series $f(x)$. After this step, the algorithm terminates. For details, see \cite[Proposition 3.2]{power}.
\end{proof}
As a filter, one can also verify a condition coming from Ehrhart theory ~\cite{class}.

\begin{remark}
    The results of \cite[Sections 4.2 and 5 ]{survey} rely on the values of the coordinates in the binomial basis, while Proposition~\ref{prop:rep} rely on the position of the non zero coordinates in the binomial basis. In that sense, both results are complementary to each other.
\end{remark}

\begin{remark}
Theorem~\ref{thm:alg}  requires a function already written as a finite linear combination of terms $\{\zetaf[n]\}_{n\in\mathbb{N}}$. The reason is that a random power series may not be expressable as a finite linear combination of the $\{\zetaf[n]\}_{n\in\mathbb{N}}$. For example, consider the function $\sum_{n=1}^\infty \sum_{j=1}^n{n\choose j} x^n$.  
\end{remark}
\begin{example}
Consider the function $8\zetaf[13]+25\zetaf[14]+18\zetaf[15]$. The algorithm of Theorem~\ref{thm:alg} efficiently tests posets until it finds a candidate. Suppose the algorithm finds $Q=D(8)*D(1)$. Then with \cite[Section 5.3]{survey} we obtain $Q_2=D(D(2))*\chain[7]$ as a Doppelg\"anger of $Q$. Up to commutativity of $*$, these are all the Doppelg\"anger Wix\'arika posets.
\end{example}

\section{Farewell}\label{Section:final}
In this note, we aimed to show how operads can relate objects of different nature such as posets and power series. 
We also gave an example in which the topological information of the operad can be used to improve an algorithm.

What other operads can interested readers study? There is an encyclopedia of operads, a project in progress by Frédéric Chapoton: Operadia \url{https://operadia.pythonanywhere.com/table}, and the Encyclopedia of types of algebras~\cite{Enc} by Jean-Louis Loday.
If one aims to study a new operad, the key point is to have a notion of composition, for example, cylinders can be composed vertically. A variant of the operad of cylinders encodes all the operations one expects in calculus, as described in~\cite[Theorem 2]{cyl}. 

The following list is not exhaustive, but a sample to illustrate papers with applications of operad theory. After all, in February 2026 the website zbMATH returns 264 papers with the word ``operad" in the title. 
Operad theory started with the work of J. M. Boardman, R. M. Vogt~\cite{Op2}, and P. May~\cite{Op1} in the field of algebraic  topology. M. Kontsevich~\cite{Op4}, D. Tamarkin~\cite{Op3}, and T. Willwacher~\cite{Op9} popularized them in mathematical physics. 
 F. Chapoton~\cite{OpTre} described the relationship between operads and combinatorial species,  S. Giraudo~\cite{Op5} studied combinatorial problems as algebras over operads.
 F. Fauvet, L. Foissy, and  D. Manchon~\cite{Op6} studied the operadic structures on finite posets. J. Baez, J. Foley, J. Moeller, and B. S. Pollard studied network models via operads ~\cite{Op7}.
 
The reader may be surprised by the diversity of fields in which operad theory is used. This is however quite natural, since operads are a categorical concept, and the aim of category theory is to study mathematics from a universal point of view. 

\section*{Acknowledgments}
An old version of this note was used to give a minicourse in Govt.\ College University Lahore in Pakistan 2022, and at the BIRS-CMO workshop ``6th meeting of the Mexican Mathematicians in the World'' 2024. The work of E. Dolores Cuenca was supported by the Korea National Research Foundation (NRF) grant funded by the Korean government (MSIT) (RS-2025-00517727).

\bibliographystyle{abbrv}
\bibliography{references}
\end{document}